\documentstyle{amsppt}
\magnification=1200
\hsize=150truemm
\vsize=224.4truemm
\hoffset=4.8truemm
\voffset=12truemm

\NoRunningHeads

\define\C{{\bold C}}
\define\R{{\bold R}}
\let\thm\proclaim
\let\fthm\endproclaim
 


\newcount\tagno
\newcount\secno
\newcount\subsecno
\newcount\stno
\global\subsecno=1
\global\tagno=0
\define\ntag{\global\advance\tagno by 1\tag{\the\tagno}}

\define\sta{\ 
{\the\secno}.\the\stno
\global\advance\stno by 1}

\define\stas{\the\stno
\global\advance\stno by 1}

\define\sect{\global\advance\secno by 1
\global\subsecno=1\global\stno=1\
{\the\secno}. }

\def\nom#1{\edef#1{{\the\secno}.\the\stno}}
\def\inom#1{\edef#1{\the\stno}}
\def\eqnom#1{\edef#1{(\the\tagno)}}

\newcount\refno
\global\refno=0

\def\nextref#1{
      \global\advance\refno by 1
      \xdef#1{\the\refno}}

\def\bref {\ref\global\advance\refno by 1\key{\the\refno}}

\nextref\CHA 
\nextref\GRO
\nextref\LEL

\topmatter

\title 
On a result by Y. Groman and J. P. Solomon
 \endtitle

\author  Julien Duval\footnote""{laboratoire de math\'ematiques, universit\'e Paris-Sud, 91405 Orsay cedex, France \newline julien.duval\@math.u-psud.fr\newline}
\footnote""{keywords : J-holomorphic curves, isoperimetric inequality \newline AMS class. : 32Q65.\newline}
\endauthor

\abstract\nofrills{\smc Abstract.} We give a short proof of a reverse isoperimetric inequality due to Y. Groman and J. P. Solomon.
\endabstract

\endtopmatter 

\document

\subhead 1. Introduction \endsubhead

\null
 
Let $(X, J)$ be a compact almost complex manifold equipped with a hermitian metric and $S\subset X$ a compact totally real submanifold of maximal dimension. Then $J$-holomorphic curves with boundary in $S$ satisfy a reverse isoperimetric inequality. Namely there exists a constant $A>0$ such that for any compact J-holomorphic curve $(C,\partial C)\subset (X,S)$ then $$\text{long}(\partial C)\leq A\ \text{area}(C).$$ This statement is due to Y. Groman and J. P. Solomon [\GRO] (with an extra term involving the genus of $C$ on the right hand side). We refer also to [\GRO] for motivation and applications. The proof of [\GRO] is geometric and combinatorial in nature. We propose here an analytic approach of this inequality based on a monotonicity principle, which gives in fact a stronger semi-local statement.

\thm{Theorem} There exists $A>0$ such that for any compact $J$-holomorphic curve $(C,\partial C)\subset(X,S)$ then 
$\text{long}(\partial C)\leq \frac A r \ \text{area}(C\cap U_r)$ where $U_r$ is the $r$-neighborhood of $S$, for $r>0$ small enough.
\fthm

Our approach is reminiscent of Lelong method [\LEL] for proving the following inequality. Let $C$ be a holomorphic curve in $\C ^n$, $m$ its multiplicity at $0$ and $B_r$ the ball of radius $r$ centered at $0$, then $$m \leq \frac {\text{area }(C\cap B_r)}{\pi r^2}.$$ It is based on the plurisubharmonicity of $\log \vert z \vert^2$. It implies that the function $a(r)=\frac 4 {r^2}\text{ area}(C\cap B_r)=\frac 1 {r^2}\int_{C\cap B_{r}}dd^c\vert z \vert^2$ is increasing. Indeed, by Stokes theorem, $a(r) =\int_{C\cap \partial B_{r}}d^c\log \vert z \vert^2$, so $a(r)-a(s)=\int_{C\cap (B_{r}\setminus B_s)}dd^c\log \vert z\vert^2\geq 0$ where $r>s>0$. As $\lim\limits_{s\to 0}a=4\pi m$ the inequality follows.

In our situation the local model is $\R ^n\subset \C ^n$ where $\R ^n=(y=0)$ with the usual coordinates $z=x+iy$ in $\C^n$. Now $\vert y\vert^2$ is strictly plurisubharmonic, so $dd^c\vert y\vert^2$ restricts to an area form on holomorphic curves. It turns out that $\vert y \vert$ is still plurisubharmonic and will play the role of $\log \vert z \vert^2$ above. This can be globalized. Near $S$ there exists a function $\rho$ looking like the square of the distance to $S$, strictly $J$-plurisubharmonic and such that $\sqrt \rho$ is still $J$-plurisubharmonic. Let us enter the details. 

\null

\subhead 2. Proof of the theorem \endsubhead

\null All objects are supposed smooth except otherwise mentioned. Recall that $S$ is a compact totally real submanifold in $X$ of maximal dimension (say $n$). This means that $TS\oplus JTS = TX\vert_S$. The point is the following  
\thm{Lemma} Near $S$ there exists a strictly $J$-plurisubharmonic function $\rho\geq0$ of class $C^2$, vanishing exactly on $S$, such that $\sqrt \rho$ is $J$-plurisubharmonic outside $S$.
\fthm

This means that $dd^J\rho>0$ and $dd^J\sqrt \rho\geq0$ where $d^J g$ stands for $-d g\circ J$. Recall that a 2-form $\theta$ is non negative (resp. strictly positive) if $\theta(v,Jv)\geq0$ (resp. $>0$) for any tangent vector $v \ne 0$.
Assuming this lemma for a while let us prove the theorem.

We may take $dd^J\rho$ as the area form of our hermitian metric near $S$. As $\rho$ is comparable to the square of the distance to $S$ we may also take $U_r=(\rho\leq r^2)$. This will only change the constant $A$ in the end.  

Take $C$ a compact J-holomorphic curve of $X$ with boundary in $S$. Precisely it is the image of a map $f:(\Sigma,i)\to (X,J)$ where $\Sigma$ is a compact Riemann surface with boundary, such that $df\circ i=J\circ df$ and $f(\partial \Sigma)\subset S$. All the integrals below should be meant parametrized by $f$, though we write them on $C$ for simplicity. 

As above $a(r)=\frac 1 r \text{ area}(C\cap U_r)=\frac 1 r \int_{C\cap U_r}dd^J\rho$ is increasing. Indeed, by Stokes theorem, $a(r)=2\int_{C\cap \partial U_r}d^J\sqrt \rho$, so $a(r)-a(s)=2\int_{C\cap(U_r\setminus U_s)}dd^J\sqrt\rho\geq0$ where $r>s>0$. Hence $\lim\limits_{s\to 0}a\leq\ a(r)$. 

On the other hand, as $\rho$ has a minimum along $S$, there exists $A>0$ such that $\vert \nabla \rho \vert \leq A s$ in $U_s$. So $A\ a(s)\geq \frac 1 {s^2} \int_{C\cap U_s}\vert \nabla (\rho\vert_C) \vert dd^J\rho = \frac 1 {s^2}\int_0^{s^2}\text{long}(C\cap (\rho=t))\ dt$ by the coarea formula (see for instance [\CHA]). Hence $A\lim\limits_{s\to 0}a\geq \text{long}(\partial C)$.
 
We conclude that $\text{long}(\partial C)\leq A\ a(r)$.

\newpage

\noindent
{\bf Proof of the lemma}. Take any function $\rho\geq 0$ near $S$, vanishing on $S$ and non degenerate transversally to $S$. It is known that $\rho$ is strictly $J$-plurisubharmonic (see below). Now $\sqrt \rho$ is not necessarily $J$-plurisubharmonic outside $S$ but we will find $B>0$ such that $\sqrt \rho + B\rho$ is. This will do for our lemma replacing $\rho$ by $(\sqrt \rho +B\rho)^2$. 

So we need only check that $dd^J\sqrt \rho \geq O(1)$. We verify it locally.

Parametrize a piece of $S$ by a piece of $\R^n$ via $\phi$. Extend $\phi$ to a local diffeomorphism from $\C^n$ to $X$ such that $d\phi\circ i=J\circ d\phi$ on $\R^n$. This amounts to prescribing the normal derivative of $\phi$ along $\R^n$. Transport the situation via $\phi$ in $\C^n$. Locally we get a function $\rho\geq 0$, vanishing on $\R^n$ and non degenerate transversally, and an almost structure $J$ coinciding with $i$ on $\R^n$. Take the usual coordinates $z=x+iy$ in $\C^n$ such that $\R^n=(y=0)$.

As $J-i=O(\vert y \vert)$ and $\rho=O(\vert y\vert^2)$ we infer that $dd^J\sqrt \rho=dd^i\sqrt\rho+ O(1)$. Note that $d^i g$ is nothing but the more familiar $d^cg$. So it is enough to check that $dd^c\sqrt\rho\geq O(1)$ where the positivity is meant with respect to $i$. Now by assumption $\rho=q+O(\vert y \vert^3)$ where $q=\sum a_{kl}(x)y_ky_l$ with $(a_{kl})$ symmetric positive definite. We may then replace $\rho$ by $q$ in our inequality. Actually we don't even have to differentiate the coefficients $a_{kl}$ of $q$ as we work modulo $O(1)$. 

So everything boils down to proving that $dd^c\sqrt q\geq 0$ where $q$ is now a {\it constant} positive definite quadratic form in $y$. By a linear change of coordinates this reduces further to the model case $q=\vert y \vert^2$. Computing we get 
$4\vert y\vert ^{3}dd^c\vert y\vert=2\vert y\vert^2dd^c\vert y\vert^2-d\vert y\vert^2\wedge d^c\vert y\vert^2=2\sum\limits_{kl}(y_kdy_l-y_ldy_k)\wedge(y_kd^cy_l-y_ld^cy_k)\geq0.$

In the same way the strict $J$-plurisubharmonicity of $\rho$ near $S$ reduces to the strict plurisubharmonicity of $\vert y \vert^2$ which is clear.

\Refs

\widestnumber\no{99}
\refno=0 
\bref \by I. Chavel \book Riemannian geometry. A modern introduction \publ Cambridge Univ. Press \yr 2006 \publaddr Cambridge
\endref 

\bref \by Y. Groman and J. P. Solomon \paper A reverse isoperimetric inequality for J-holomorphic curves \jour Geom. Funct. Anal. \vol24\yr2014\pages1448-1515
\endref

\bref \by P. Lelong \paper Propri\'et\'es m\'etriques des vari\'et\'es analytiques d\'efinies par une \'equation \jour Ann. Sci. Ecole Norm. Sup. \vol67\yr1950\pages393--419
\endref

\endRefs

\enddocument